\documentclass[11pt]{article}

\usepackage{mathrsfs}
\usepackage{latexsym}
\usepackage{amsmath,amsthm}
\usepackage{amsfonts}
\usepackage{url}

\usepackage{amssymb}
\usepackage{pifont}

\newcommand{\h}{\eta}
\newcommand{\Vbar}{\bar{V}}

\usepackage{bbm}
\newcommand{\chr}{\boldsymbol{\mathbbm{1}}} 
\newcommand{\pred}[1]{\chr_{\left\{ #1 \right\}}}
\newcommand{\prs}{\vec{P}}
\newcommand{\pr}[1]{\prs\!\tlprn{#1}}
\newcommand{\mexp}{\vec{E}}

\newcommand{\E}{\mexp}
\renewcommand{\P}{\prs}

\newcommand{\Dn}{\Delta_n}

\newcommand{\TV}[1]{\nrm{#1}_{\textrm{{\tiny \textup{TV}}}}}

\newcommand{\Lipw}[2]{\nrm{#2}_{\textrm{{\tiny \textup{Lip}}},#1}}

\newcommand{\diam}{\operatorname{diam}}

\newcommand{\subplus}{_{
+
}}

\newcommand{\pl}[1]{\paren{#1}\subplus}

\usepackage{geometry}
\makeatletter
\theoremstyle{plain}
\newtheorem{thm}{Theorem}[section]
\theoremstyle{plain}    
\newtheorem*{thm*}{Theorem}
 \theoremstyle{remark}
 \newtheorem{rem}[thm]{Remark}
 \theoremstyle{plain}    
 \newtheorem{lem}[thm]{Lemma} 
 \theoremstyle{plain}    
 \newtheorem{cor}[thm]{Corollary} 
\makeatother
\newcommand{\bethn}{\begin{thm}}
\newcommand{\enthn}{\end{thm}}
\renewcommand{\beth}{\begin{thm*}}
\newcommand{\enth}{\end{thm*}}
\newcommand{\bepf}{\begin{proof}}
\newcommand{\enpf}{\end{proof}}
\newcommand{\belen}{\begin{lem}}
\newcommand{\enlen}{\end{lem}}
\newcommand{\becon}{\begin{cor}}
\newcommand{\encon}{\end{cor}}
\newcommand{\ben}{\begin{enumerate}}
\newcommand{\een}{\end{enumerate}}
\newcommand{\bit}{\begin{itemize}}
\newcommand{\eit}{\end{itemize}}

\renewcommand{\vec}[1]{\bs{\mathrm{#1}}}

\newcommand{\phiwnorm}[2]{\nrm{#1}_{\Phi,#2}}
\newcommand{\psiwnorm}[2]{\nrm{#1}_{\Psi,#2}}
\newcommand{\basicspace}{{\cal S}}
\newcommand{\X}{\basicspace}
\newcommand{\supr}[1]{^{(#1)}}
\newcommand{\supl}[1]{^{+#1}}

\newcommand{\sseq}[3]{#1_{#2}^{#3}}  

\newcommand{\sd}{}
\newcommand{\scat}[1]{#1}

\newcommand{\nrm}[1]{\left\Vert #1 \right\Vert}

\newcommand{\iprod}[2]{\left\langle #1 , #2 \right\rangle}

\newcommand{\f}{\varphi}
\renewcommand{\k}{\kappa}

\newcommand{\psin}[1]{\Psi_{#1}}
\newcommand{\psiwn}[2]{\Psi_{{#1},{#2}}}
\newcommand{\phiwn}[2]{\Phi_{{#1},{#2}}}

\newcommand{\calL}{\mathcal{L}}

\newcommand{\calX}{\mathcal{X}}

\newcommand{\RR}{\mathbb{R}}
\newcommand{\NN}{\mathbb{N}}

\newcommand{\beq}{\begin{eqnarray*}}
\newcommand{\eeq}{\end{eqnarray*}}
\newcommand{\beqn}{\begin{eqnarray}}
\newcommand{\eeqn}{\end{eqnarray}}
\newcommand{\paren}[1]{\left( #1 \right)}
\newcommand{\sqprn}[1]{\left[ #1 \right]}
\newcommand{\tlprn}[1]{\left\{ #1 \right\}}
\newcommand{\set}[1]{\tlprn{#1}}
\newcommand{\abs}[1]{\left| #1 \right|}

\newcommand{\gn}{\, | \,}

\newcommand{\ts}{\textstyle}
\newcommand{\bs}{\boldsymbol}

\newcommand{\hide}[1]{}
\newcommand{\oo}[1]{\frac{1}{#1}}

\newcommand{\defeq}{\doteq}

\title{A Linear Programming Inequality
with Applications to Concentration of Measure}
\author{Leonid Kontorovich\\
School of Computer Science\\ 
Carnegie Mellon University\\ 
Pittsburgh, PA 15213\\
USA \\
\url{lkontor@cs.cmu.edu}
}
\begin{document}
\maketitle
\abstract{
We prove 
an elementary
yet useful inequality bounding the maximal value of
certain linear programs. 
This leads directly
to a bound on the martingale
difference for 
arbitrarily
dependent random variables, 
providing
a generalization
of some recent concentration of measure results. The linear
programming inequality may be of independent interest.
}

\newcommand{\rp}{\RR_+}
\newcommand{\rpp}{\RR_{++}}
\section{Introduction}
\subsection{Background}
Over the past decade
there
has been a 
flurry of new concentration of measure
inequalities; we refer the reader to
\cite{ledoux01} for an in-depth survey, or
\cite{kontram06,kont06-metric-mix,lugosi} for some more recent advances. 

In \cite{kontram06} the martingale difference method was employed in a
novel way to obtain a general concentration inequality for dependent
random variables, with respect to the (unweighted) Hamming metric. At
the core of that approach lies a certain linear programming inequality
associated with bounding martingale differences 
\cite[Theorem 4.8]{kontram06}. In this paper, we give a considerably
simpler proof of 
a rather more general
result, 
extending it
to 
the
weighted Hamming
metrics. The applications to measure concentration are immediate 
(culminating in Corollary~\ref{cor:concmeas});
additionally, it is hoped that the linear
programming inequality and the 
technique
employed for proving it
will find
further applications.

Since the main focus of this paper is the inequality in
Theorem~\ref{thm:phipsi}, we forgo a detailed discussion of measure
concentration and how our bound relates to existing results. Such a
discussion may be found in 
\cite{kontram06,kont06-metric-mix}.

\subsection{Notational conventions}
\label{sec:notconv}
Throughout this paper, 
$\X$ will
denote a finite set.
Random variables are capitalized ($X$), specified sequences (words)
are written in lowercase ($x\in\X^n$), the shorthand
$\sseq{X}{i}{j}\defeq
(X_i,\ldots,X_j)$ is used for all sequences, and
word concatenation is denoted using the multiplicative notation:
$\sseq{x}{i}{j} \sseq{x}{j+1}{k}=\sseq{x}{i}{k}$.
Similarly, if $w\in\RR^n$ and $1\leq k\leq\ell\leq n$, then
$\sseq{w}{k}{\ell}\defeq
(w_k,\ldots,w_\ell)
\in\RR^{k-\ell+1}$.

We use the indicator variable
$\pred{\cdot}$ 
to assign 0-1 truth values 
to the predicate in 
$\set{\cdot}$.
The ramp function
is defined by $\pl{z}=z\pred{z>0}$.
The 
positive
reals are denoted by $\rp\defeq(0,\infty)$.

The probability 
$\P$
and expectation $\E$ operators
are defined 
with respect 
the measure space specified in context.

\section{Linear programming inequality}
We begin with a natural generalization of some of the definitions in
\cite{kontram06}. Fix a finite set $\X$,
$n\in\NN$ and $w\in\rp^n$. Then
\ben
\item $K_n$ denotes the set of all functions $\k:\X^n\to\RR$
(and
$K_0\doteq\RR$)
\item the {\it weighted Hamming metric} on $\X^n\times\X^n$ is defined
      by
\beqn
\label{eq:wham}
d_w(x,y) &=& \sum_{i=1}^n w_i\pred{x_i\neq y_i}
\eeqn
\item for $\f\in K_n$, its {\em Lipschitz constant} with respect to
  $d_w$, denoted by $\Lipw{w}{\f}$, is defined to be the smallest $c$
  for which
\beq
\abs{\f(x)-\f(y)} &\leq& c d_w(x,y)
\eeq
for all $x,y\in\X^n$; any $\f$ with $\Lipw{w}{\f}\leq c$ is called $c$-Lipschitz
\item 
for $v\in[0,\infty)$,
define $\phiwn{w}{n}\supl v\subset K_n$ to be the set of all\hide{
\beq
\abs{\f(x)-\f(y)} &\leq& d_w(x,y)
\eeq
for all $x,y\in\X^n$;}
$\f$ such that 
$\Lipw{w}{\f}\leq1$ and
\beq
0\;\leq\; \f(x) \;\leq\; v+\sum_{i=1}^n w_i,
\qquad x\in\X^n;
\eeq
we omit the $+v$ superscript when $v=0$, writing
simply $\phiwn{w}{n}$
\item
the
{\em marginal projection operator}
$(\cdot)'$
takes $\k\in K_n$ to $\k'\in K_{n-1}$
by
\beq
\k'(y) &\doteq& \sum_{x_1\in\X} \k (x_1  y),
\qquad x \in \X^{n-1};
\eeq
for $n =1$,  $\k'$ is the scalar $\k' = \sum_{x_1 \in \X} \k (x_1)$
\item
for $y\in\X$,
the {\em $y$-section} operator $(\cdot)_y
$ takes $\k\in K_n$ to $\k_y\in K_{n-1}$ by
\beq
\k_y(x) &\doteq& \k(x  y) ,
\qquad x \in \X^{n-1};
\eeq
for $n = 1$,
$\k_y(\cdot)$ is the scalar $\k(y)$
\item
the functional $\psiwn{w}{n}:K_n\to\RR$ is defined by
$\psiwn{w}{0}(\cdot)=0$ and
\beqn
\label{eq:psiwn}
\psiwn{w}{n}(\k) &\defeq& w_1\!\sum_{x\in\X^n}\pl{\k(x)}+
\psiwn{
\sseq{w}{2}{n}
}{n-1}(\k');
\eeqn
when $w_i\equiv1$ we omit it from the subscript, writing simply
$\psin{n}$
\item 
the finite-dimensional vector space
$K_n$ is equipped 
with the inner product
\beq
\iprod{\k}{\lambda} &\doteq& \sum_{x\in\X^n}\k(x)\lambda(x)
\eeq

\item
two norms are defined on $\k\in K_n$:
the $\Phi_w$-norm,
\beqn
\label{eq:phinorm}
\phiwnorm{\k}{w} &\defeq&
\sup_{\f\in\phiwn{w}{n}}\abs{\iprod{\k}{\f}}
\eeqn
and the $\Psi_w$-norm,
\beqn
\label{eq:phinorm}
\psiwnorm{\k}{w} &\defeq&
\max_{s=\pm1} \psiwn{w}{n}(s\k).
\eeqn

\een
\begin{rem}
For the special case
$w_i\equiv 1$,
$d_w$ is the unweighted Hamming metric used in
\cite{kontram06}.
It is straightforward to verify that $\Phi_w$-norm and
$\Psi_w$-norm satisfy the vector-space norm axioms
for any $w\in\rp^n$; this is done in 
\cite{kontram06} for $w_i\equiv 1$. Since we will not be appealing to any
norm properties of these functionals, we omit the proof.
Note that for any $y \in \X$, the marginal projection and $y$-section
operators commute; in other words,  for
$\k\in K_{n+2}$, we have $(\k')_y=(\k_y)' \in K_n$ and so
we can denote this common value 
by $\k_y'\in K_n$: 
\beq
\k_y'(z) = \sum_{x_1\in\X} \k_y(x_1  z)
         = \sum_{x_1\in\X} \k  (x_1  z  y),
\qquad z\in\X^n
. 
\eeq
\end{rem}

\hide{
\begin{verbatim}
-reference to LP is incidental as no results from LP are used
-for lack of a better name, we have referred to thm:main as the ``LP
-ineq'' and the ``phi-psi'' ineq.

-S is finite set
-projection, y-section
-\sseq notation for w
\end{verbatim}

All undefined notation is found in 
\cite{kontram06}. We will use
$\rp\defeq(0,\infty)$.
Let $\X$ be a finite set.

For $\k\in K_n$, let $\k'\in K_{n-1}$ be the marginal projection of
$\k$. 
}

The main result of this section is
\bethn
\label{thm:main}
For 
all 
$w\in\rp^n$
and
all $\k\in K_n$, we have
\beqn
\label{eq:phipsinorm}
\phiwnorm{\k}{w} &\leq& \psiwnorm{\k}{w}.
\eeqn
\enthn
\begin{rem}
We refer to (\ref{eq:phipsinorm}) -- more properly, to 
(\ref{eq:mainclaim}), from which the former immediately follows -- as
a
{\em linear programming inequality} for the reason that 
$F(\cdot)=\iprod{\k}{\cdot}$ is a linear function being maximized
over the finitely generated, compact, convex polytope
$\phiwn{w}{n}\subset\RR^{\X^n}$.
We make no use of this
simple fact and therefore forgo its proof, but see
\cite[Lemma 4.4]{kontram06}
for a proof of a closely related claim.
The term ``linear programming'' is a bit of a red herring since no
actual LP
techniques are being used; for lack of an obvious natural 
name, we have alternatively referred
to 
(\ref{eq:phipsinorm})
in previous papers and talks as
the ``$\Phi$-norm bound'' or the ``$\Phi$-$\Psi$ inequality.''
\end{rem}

The key technical lemma is 
a
decomposition of 
$\psiwn{w}{n}(\cdot)$ in terms of $y$-sections, proved in
\cite{kontram06} for the case $w_i\equiv 1$:

\belen
\label{lem:psidecomp}
For all $n\geq1$,
$w\in\rp^{n}$
and $\k\in K_{n}$,
we have
\beqn
\label{eq:psidecomp}
\psiwn{w}{n}(\k) & = &  \sum_{y\in\X} \sqprn{
\psiwn{
\sseq{w}{1}{n-1}
}{n-1}(\k_y) +
w_{n}\!\pl{ \sum\limits_{x\in\X^{n-1}}  \k_y(x) }}.
\eeqn
\enlen
\hide{
\begin{rem}
The claim holds as stated for $n=1$, but we must recall that $\X^0$ is
the set containing a single (null) word 
and that for
$\k\in K_1$, $\k_y\in K_0$ is the scalar $\k(y)$.
\end{rem}
}
\bepf
We proceed by induction on $n$. 
To prove the $n=1$ case,
recall that $\X^0$ is
the set containing a single (null) word 
and that for
$\k\in K_1$, $\k_y\in K_0$ is the scalar $\k(y)$. Thus, by definition
of 
$\psiwn{w}{1}(\cdot)$, we have
\beq
\psiwn{w}{1}(\k) & = &  w_1\!\sum_{y\in\X} \k(y),
\eeq
which proves (\ref{eq:psidecomp}) for $n=1$.

Suppose the claim holds for some $n=\ell\geq1$. Pick any
$w\in\rp^{\ell+1}$ and $\k\in K_{\ell+1}$ and examine

\beqn
\sum_{y\in\X} \sqprn{
\psiwn{\sseq{w}{1}{\ell}}{\ell}(\k_y) + 
w_{\ell+1}\!\pl{\sum_{x\in\X^{\ell}} \k_y(x)}
}
&&\nonumber\\
&&\nonumber
\hspace{-3cm}
=\;
\sum_{y\in\X}\sqprn{
\paren{
w_1\!\sum_{x\in\X^{\ell}} \pl{\k_y(x)} 
+
\psiwn{\sseq{w}{2}{\ell}}{\ell-1}(\k'_y)} 
+ 
w_{\ell+1}\!\pl{\sum_{x\in\X^{\ell}} \k_y(x)}
}
\\\nonumber\\
&&
\hspace{-3cm}
=\;
\label{eq:psineqrhs1}
\sum_{y\in\X}\sqprn{
\psiwn{\sseq{w}{2}{\ell}}{\ell-1}(\k'_y) +
w_{\ell+1}\!\pl{\sum_{u\in\X^{\ell-1}} \k'_y(u)}
}
+ 
w_1\!\sum_{z\in\X^{\ell+1}} \pl{\k(z)}
\nonumber\\
\eeqn
where the first equality follows from the definition of
$\psiwn{\sseq{w}{1}{\ell}}{\ell}$ in (\ref{eq:psiwn}) and the second one
from the easy 
identities
\beq
\sum_{y\in\X}\sum_{x\in\X^{\ell}} \pl{\k_y(x)}
&=&
\sum_{z\in\X^{\ell+1}} \pl{\k(z)}
\eeq
and
\beq
\sum_{x\in\X^{\ell}} \k_y(x) &=& \sum_{u\in\X^{\ell-1}} \k'_y(u).
\eeq

On the other hand, by definition we have
\beqn
\label{eq:psineqrhs2}
\psiwn{w}{\ell+1}(\k) 
&=&
w_1\sum_{z\in\X^{\ell+1}}\pl{\k(z)} + \psiwn{\sseq{w}{2}{\ell+1}}{\ell}(\k').
\eeqn
To compare 
the r.h.s. of (\ref{eq:psineqrhs1}) with
the r.h.s. of (\ref{eq:psineqrhs2}), note that the
$w_1\sum_{z\in\X^{\ell+1}}\pl{\k(z)}$
term 
is common to both
and
\beq
\sum_{y\in\X}\sqprn{
\psiwn{\sseq{w}{2}{\ell}}{\ell-1}(\k'_y)
+ 
w_{\ell+1}\pl{\sum_{u\in\X^{\ell-1}} \k'_y(u)}
}
&=&
\psiwn{\sseq{w}{2}{\ell+1}}{\ell}(\k')
\eeq
by the inductive hypothesis. This establishes 
(\ref{eq:psidecomp}) for $n=\ell+1$ and proves the claim.
\enpf

Our main result, Theorem~\ref{thm:main}, is an immediate consequence
of

\bethn
\label{thm:phipsi}
For all $n\geq1$,
$w\in\rp^{n}$,
$v\in[0,\infty)$
and $\k\in K_{n}$,
we have
\beqn
\label{eq:mainclaim}
\sup_{\f\in\phiwn{w}{n}\supl v} {\iprod{\k}{\f}}
&\leq & \psiwn{w}{n}(\k)
+ v\pl{\sum_{x\in\X^n}\k(x)}
.
\eeqn
\enthn
\bepf
We will prove the claim by induction on $n$. For $n=1$, pick any
$w_1\in\rp$, 
$v\in[0,\infty)$ and
$\k\in K_1$. Since by construction any
$\f\in \phiwn{w_1}{1}\supl v$ is $w_1$-Lipschitz with respect to the
discrete
metric on $\X$, $\f$ must be of the form
\beq
\f(x) &=& \tilde\f(x) + \tilde v,
\qquad x\in\X,
\eeq
where $\tilde\f:\X\to[0,w_1]$ and $0\leq\tilde v\leq v$ 
(in fact, we have the explicit value
$\tilde v=\pl{\max_{x\in\X}\f(x)-w_1}$).
Therefore,
\beqn
\label{eq:kptil}
\iprod{\k}{\f} &=& \iprod{\k}{\tilde\f} + \tilde v\sum_{x\in\X}\k(x).
\eeqn
The first term 
in the r.h.s.
of (\ref{eq:kptil}) is clearly maximized when
$\tilde\f(x)=w_1\pred{\k(x)>0}$ for all $x\in\X$, which shows that it
is bounded by $\psiwn{w_1}{1}(\k)$. Since the second term 
in the r.h.s.
of
(\ref{eq:kptil}) is bounded by $v\pl{\sum_{x\in\X}\k(x)}$, we have established
(\ref{eq:mainclaim}) for $n=1$.

Now suppose the claim holds
for $n=\ell$, and pick any
$w\in\rp^{\ell+1}$, 
$v\in[0,\infty)$ and
$\k\in K_{\ell+1}$.
By the reasoning given above
(i.e., using the fact that $0\leq\f\leq v+\sum_{i=1}^{\ell+1}w_i$
and 
that
$\f$ is $1$-Lipschitz with respect to $d_w$),
any 
$\f\in\phiwn{w}{\ell+1}\supl v$, 
must be of the form $\f=\tilde\f + \tilde v$, where
$\tilde\f\in\phiwn{w}{\ell+1}$ and $0\leq\tilde v\leq v$.
Thus we write 
$\iprod{\k}{\f}=\iprod{\k}{\tilde\f}+
\tilde v\sum_{x\in\X^{\ell+1}}\k(x)$
and
decompose
\beqn
\label{eq:kfy}
\iprod{\k}{\tilde\f} &=& \sum_{y\in\X} \iprod{\k_y}{\tilde\f_y},
\eeqn
making
the obvious but 
crucial
observation that
\beq
\tilde\f\in \phiwn{w}{\ell+1}
&\Longrightarrow&
\tilde\f_y
\in \phiwn{\sseq{w}{1}{\ell}}{\ell}\supl{w_{\ell+1}}.
\eeq
Then it follows by the inductive hypothesis that
\beqn
\label{eq:mainind}
\iprod{\k_y}{\tilde\f_y} &\leq&
\psiwn{\sseq{w}{1}{\ell}}{\ell}(\k_y)
+ w_{\ell+1}\!\pl{\sum_{x\in\X^\ell}\k_y(x)}
.
\eeqn
Applying Lemma~\ref{lem:psidecomp}
to 
(\ref{eq:mainind}), we have
\beqn
\label{eq:psiy}
\sum_{y\in\X} \iprod{\k_y}{\tilde\f_y}
\;\leq\;
\sum_{y\in\X} \sqprn{
\psiwn{\sseq{w}{1}{\ell}}{\ell}(\k_y)
+ w_{\ell+1}\!\pl{\sum_{x\in\X^\ell}\k_y(x)}
}
\;=\;
\psiwn{w}{\ell+1}(\k).
\eeqn
This, 
combined with
(\ref{eq:kfy}) and 
the trivial bound
\beq
\tilde v\sum_{x\in\X^{\ell+1}}\k(x)
&\leq&
v\pl{\sum_{x\in\X^{\ell+1}}\k(x)}
\eeq
proves the
claim for $n=\ell+1$ and 
hence
for all $n$.
\enpf

\section{Applications to concentration of measure}
\label{sec:apps}
This section assumes some familiarity with the notion of measure
concentration; see the References section (in particular,
\cite{ledoux01,lugosi}) for introductory and survey material.
Briefly, we shall concern ourselves with the metric probability space
$(\X^n,d_w,\P)$ where $\X$ is a finite set, $w\in\rp^n$, $d_w$ is the
weighted Hamming metric defined in (\ref{eq:wham}) and $\P$ is a
(possibly non-product) probability measure on $\X^n$. For random
variables $f:\X^n\to\RR$, our goal is to bound
$\pr{\abs{f-\E f}>t}$.

The method of martingale differences has been used to prove
concentration of measure results since the work of Hoeffding, Azuma,
and McDiarmid; see the exposition and references in
\cite{kontram06,kont06-metric-mix}.
Let $(\X^n,d_w,\P)$ be as defined above and associate to it the
(canonical) random
process 
$(X_i)_{1\leq i\leq n}$, $X_i\in\X$, satisfying
\beq
\pr{X\in A} &=& \P(A)
\eeq
for any $A\subset\X^n$. 

For 
$1\leq i\leq n$,
$f:\X^n\to\RR$
and $\sseq{y}{1}{i}\in\X^i$, 
define the
{\em martingale difference}
\beqn
\label{eq:videf}
V_i(f;\sseq{y}{1}{i}) &=&
\E[f(X)\gn\sseq{X}{1}{i}=\sseq{y}{1}{i}] -
\E[f(X)\gn\sseq{X}{1}{i-1}=\sseq{y}{1}{i-1}].
\eeqn
Let
\beqn
\label{eq:Vbar}
\Vbar_i(f) &\defeq& \max_{\sseq{y}{1}{i}\in\X^i} \abs{V_i(f;\sseq{y}{1}{i})}
\eeqn
and
\beq
D^2(f) &\defeq& \sum_{i=1}^n \Vbar_i^2(f).
\eeq
Then Azuma's inequality~\cite{azuma} states that
\beqn
\label{eq:azuma}
\pr{\abs{f-\E f}>t} &\leq& 2\exp(-t^2/2D^2(f))
\eeqn
(see~\cite{ledoux01} for a modern presentation and a short proof of
(\ref{eq:azuma})). 

In \cite{kontram06} and \cite{kont06-metric-mix}, a technique was
developed for bounding the martingale difference $V_i(f;y)$ in terms
of the Lipschitz constant of $f$ and mixing properties of the measure
$\P$. To this end, we introduce the so-called $\eta$-mixing
coefficients (see discussion ibid. regarding the appearance of these
coefficients in earlier work of Marton~\cite{marton03} and
Samson~\cite{samson00}). 

For $1\leq i<j\leq n$
and $x\in\X^i$,
let
$$\calL(\sseq{X}{j}{n}\gn \sseq{X}{1}{i}=x)
$$ be the law 
(distribution)
of $\sseq{X}{j}{n}$ conditioned on 
$\sseq{X}{1}{i}=x$. 
For $y\in\X^{i-1}$ and 
$z,z'\in\X$, 
define
\beqn
\label{eq:hdef}
\h_{ij}(y,z,z') &=&
\TV{
\calL(\sseq{X}{j}{n}\gn \sseq{X}{1}{i}=\scat{y\sd z})-
\calL(\sseq{X}{j}{n}\gn \sseq{X}{1}{i}=\scat{y\sd z'})
},
\eeqn
where $\TV{\cdot}$ is the total variation norm, defined here, for a
signed measure $\tau$ on a finite space $\calX$ by
\beq
\TV{\tau} &\defeq& {\ts\oo2}\sum_{x\in\calX}\abs{\tau(x)}.
\eeq
Additionally, define
\beq
\bar\h_{ij} &=&
\max_{y\in\X^{i-1}} 
\max_{z,z'\in\X} 
\h_{ij}(y,z,z').
\eeq

The main application of Theorem~\ref{thm:phipsi} to measure
concentration is the following bound on the martingale difference:
\bethn
Let $\X$ be a finite set, and
let $(X_i)_{1\leq i\leq n}$, $X_i\in\X$ be the random process
associated with the 
measure $\P$ on $\X^n$.
Let $\Dn$ be the upper-triangular $n\times n$ matrix defined by
$(\Dn)_{ii}=1$ and
\beqn
\label{eq:Ddef}
(\Dn)_{ij} = \bar\h_{ij}
\eeqn
for $1\leq i<j\leq n$.
Then, for all 
$w\in\rp^n$ and
$f:\X^n\to\RR$, we have
\beqn
\label{eq:sumVi}
\sum_{i=1}^n \Vbar_i^2(f) &\leq& \Lipw{w}{f}^2\nrm{\Dn w}_2^2
\eeqn
where $\Vbar_i^2(f)$ is defined in (\ref{eq:Vbar}).
\enthn
\begin{rem}
Since
$\Vbar_i(f)$
and
$\Lipw{w}{f}$
are both homogeneous functionals of $f$ (in the sense that
$T(af)=|a|T(f)$ for $a\in\RR$), there is no loss of generality in
taking $\Lipw{w}{f}=1$. 
Additionally, since
$V_i(f;y)$ is
translation-invariant (in the sense that $V_i(f;y)=V_i(f+a;y)$ for
all $a\in\RR$), there is no loss of generality in restricting the
range of $f$ to $[0,\diam_{d_w}(\X^n)]$. In other words, it suffices
to consider $f\in\phiwn{w}{n}$.
Since essentially this result (for $w_i\equiv1$) is proved in
\cite{kontram06} in some detail, we only give a proof sketch here,
highlighting the changes needed for general $w$. We also remark that
the extension of this result to countable $\X$ is quite
straightforward, along the lines of~\cite[Lemma 6.1]{kontram06}.
\end{rem}
\bepf
It was shown in Section 5 of~\cite{kontram06} that if $d_w$ is the unweighted
Hamming metric (that is, $w_i\equiv 1$) and $f:\X^n\to\RR$ is
$1$-Lipschitz with respect to $d_w$, then
\beqn
\label{eq:Vhbound}
\Vbar_i(f) &\leq& 1+ \sum_{j=i+1}^n \bar\h_{ij}.
\eeqn
This was done by showing that for $1\leq i\leq n$ and $y\in\X^i$,
there is a $g_i:\X^n\to\RR$ (whose 
explicit
construction,
depending on 
$y$ and
$\P$,
is given 
\cite[Eq. (5.2)]{kontram06}), such that for all $f:\X^n\to\RR$, we have
\beqn
\label{eq:Viprodbound}
\abs{V_i(f;y)} &\leq& \abs{\iprod{g_i}{f}}.
\eeqn
It was additionally shown 
in the course of proving~\cite[Theorem 5.1]{kontram06} that
\beq
\iprod{g_i}{f} &=& \iprod{T_y g_i}{T_y f},
\eeq
where
the operator $T_y:K_n\to K_{n-i+1}$
is
defined by
\beq
(T_y h)(x) &\defeq& h(yx),
\qquad\text{for all } x\in\X^{n-i+1}.
\eeq
Appealing to
\cite[Theorem 4.8]{kontram06} -- the $w_i\equiv 1$ special
case of Theorem
\ref
{thm:phipsi} 
proved here -- 
we get
\beqn
\iprod{T_y g_i}{T_y f} &\leq & \psin{n}(T_y g_i)
.
\eeqn
It is shown in
\cite[Theorem 5.1]{kontram06} that
the form of $g_i$ implies that
\beqn
\psin{n}(T_y g_i) &\leq& 1+ \sum_{j=i+1}^n \bar\h_{ij},
\eeqn
establishing (\ref{eq:Vhbound}).\hide{
More precisely, 
let
the operator $T_y:K_n\to K_{n-i+1}$
be
defined by
\beq
(T_y h)(x) &\defeq& h(yx),
\qquad\text{for all } x\in\X^{n-i+1}.
\eeq}
To generalize (\ref{eq:Vhbound}) to $w_i\not\equiv1$, we use the fact
that
if $f\in K_n$ is $1$-Lipschitz with respect to $d_w$, then
$T_y f\in K_{n-i+1}$ is $1$-Lipschitz with respect to
$d_{\sseq{w}{i}{n}}$.
Thus,
applying Theorem~\ref{thm:phipsi},
we get
\beqn
\iprod{T_y g_i}{f} &\leq&
\psiwn{\sseq{w}{i}{n}}{n-i+1}(T_yg_i).
\eeqn
It follows directly from the definition of $\psiwn{w}{n}$ and
the calculation in \cite[Theorem 5.1]{kontram06} that
\beqn
\Vbar_i(f) &\leq& w_i+ \sum_{j=i+1}^n w_j\bar\h_{ij}\\
&=& \sum_{j=1}^n (\Dn)_{ij}w_j
= (\Dn w)_i
.
\eeqn
Squaring and summing over $i$, we obtain (\ref{eq:sumVi}).
\enpf

\becon
\label{cor:concmeas}
Let $\X$ be a 
finite
set and $\P$ a measure on $\X^n$, for
$n\geq1$. 
For any
$w\in\rp^n$ and $f:\X^n\to\RR$,
we have
\beq
\pr{\abs{f-\E f}>t} &\leq&
2\exp\paren{-\frac{t^2}{2\Lipw{w}{f}^2\nrm{w}_2^2\nrm{\Dn}_2^2}}
\eeq
where $\nrm{\Dn}_2$ is the $\ell_2$ operator norm of the matrix
defined in (\ref{eq:Ddef}).
\encon
\bepf
Since by definition of the $\ell_2$ operator norm, 
$\nrm{\Dn w}_2 \leq \nrm{\Dn}_2\nrm{w}_2$, the claim follows
immediately
via (\ref{eq:azuma}) and (\ref{eq:sumVi}).
\hide{
By definition,
\beq
\nrm{\Dn}_2 = \sup_{0\neq z\in\RR^n: \nrm{z}_2\leq 1}\nrm{\Dn z}_2
\geq \nrm{\Dn w}_2.
\eeq
Combining (\ref{eq:azuma}) and (\ref{eq:sumVi}) proves the claim.
}
\enpf

\section*{Acknowledgements}
I thank 
John Lafferty and
Kavita Ramanan for helpful discussions.


\begin{thebibliography}{10}

\bibitem{azuma}
Kazuoki Azuma, ``Weighted sums of certain dependent random variables.''
{\it Tohoku Math. Journal}, 19:357--367, 1967. 

\bibitem{kontram06}
Leonid Kontorovich and Kavita Ramanan,
``Concentration Inequalities for Dependent Random Variables via the Martingale Method.''
\url{http://arxiv.org/abs/math.PR/0609835}, 2006.

\bibitem{kont06-metric-mix}
Leonid Kontorovich, ``Metric and Mixing
Sufficient Conditions
for
Concentration
of
Measure.'' 
\url{http://arxiv.org/abs/math.PR/0610427}, 2006.

\bibitem{ledoux01}
Michel Ledoux, {\it The Concentration of Measure Phenomenon}, 
Mathematical Surveys and Monographs Vol. 89,
American
Mathematical Society, 2001.

\bibitem{lugosi}
G\'abor Lugosi, ``Concentration-of-measure inequalities.''
\url{http://www.econ.upf.es/~lugosi/anu.ps}

\bibitem{marton03}
Katalin Marton, ``Measure concentration and strong mixing.''
{\it Studia Scientiarum Mathematicarum Hungarica}, Volume 40,
Numbers 1-2, pp. 95--113(19), 2003.

\bibitem{samson00}
Paul-Marie Samson,
``Concentration of measure inequalities for Markov chains and $\Phi$-mixing processes.''
 {\it Ann. Probab.}, Vol. 28, No. 1, 416--461, 2000.


\end{thebibliography}
\end{document}